\documentclass[11pt]{article}
\usepackage{amsmath, amssymb, amscd, epsfig}

\newtheorem{thm}{Theorem}
\newtheorem{prop}[thm]{Proposition}
\newtheorem{lem}[thm]{Lemma}

\newtheorem{cor}[thm]{Corollary}
\newtheorem{rem}[thm]{Remark}
\newtheorem{df}[thm]{Definition}
\newtheorem{ex}[thm]{Example}

\renewcommand{\epsilon}{\varepsilon}

\newcommand{\BB}{\mathbb}
\newcommand{\g}{\mathfrak}
\newcommand{\pf}{\noindent {\it Proof. }}
\newcommand{\qed}{\nopagebreak $\qquad$ $\square$ \vskip5pt}
\newcommand{\separate}{\vskip5pt}
\newcommand{\VF}{\operatorname{VF}}

\newcommand{\supp}{\operatorname{supp}}

\newcommand{\re}{\operatorname{Re}}
\newcommand{\tr}{\operatorname{Tr}}

\begin{document}

\title{\bf On equivariant cohomology, integrals of equivariant forms and
Duistermaat-Heckman measures for non-compact group actions}
\author{Matvei Libine}
\maketitle

\begin{abstract}
This is an expanded version of the talk I gave at the Oberwolfach workshop
on Cohomological Aspects of Hamiltonian Group Actions and Toric Varieties,
MIT and several other universities.
For a real Lie group $G$ acting on a manifold $M$, I give definitions of
$G$-equivariant cohomology and $G$-equivariant forms on $M$.
Let $\g g$ be the Lie algebra of $G$, and let
$\alpha(X)$ be an equivariantly closed form on $M$ depending
on $X \in \g g$.
For $X \in \g g$, we denote by $M^X$ the set of zeroes of the vector
field on $M$ induced by the infinitesimal action of $X$.
Then the integral localization formula says that the integral of
$\alpha(X)$ can be expressed as a sum over the set of zeroes $M^X$
of certain {\em local} quantities of $M$ and $\alpha$:
$$
\int_M \alpha(X) = \sum_{p \in M^X}
\text{local invariant of $M$ and $\alpha$ at $p$}.
$$
For compact groups this result was proved by N.~Berline and M.~Vergne
\cite{BV1} and independently by M.~Atiyah and R.~Bott \cite{AB}
more than twenty years ago, but practically no progress had been
made until very recently in generalizing it to non-compact group actions.

I use an interplay between recent results from representation theory
and algebraic geometry to find such a generalization (\ref{mainequation}).
This generalization provides, for instance, a geometric proof of the
integral character formula from representation theory.
These results strongly suggest that many theorems which were previously known
in the compact group setting only can be generalized to non-compact groups.
The main purpose of these notes is to explain the new localization formula
and to describe the setting suitable for studying non-compact group actions.

\end{abstract}


\tableofcontents

\begin{section}
{Equivariant cohomology}
\end{section}

Let $G$ be a real Lie group acting continuously on a topological space $X$.
We make absolutely no assumption on compactness of $G$ and
we copy the definition of equivariant cohomology given in Section 1.1 of
\cite{GS} for compact groups.

\begin{df}  \label{def}
We define the {\em $G$-equivariant cohomology of $X$} to be the
ordinary cohomology of the quotient space $(X \times E)/G$:
$$
H^*_G(X) \quad =_{def} \quad H^*((X \times E)/G),
$$
where $E$ is any contractible space on which $G$ acts freely and the projection
$E \twoheadrightarrow E/G$ is a (locally trivial) fiber bundle.
\end{df}

\begin{rem}
The argument given in Section 1.1 of \cite{GS} shows that the above definition
is independent of the choice of $E$, provided that such a space exists.
If $G$ acts on $X$ freely, then $H^*_G(X) = H(X/G)$.
\end{rem}

We still need to show that there exists a contractible space $E$ on which $G$
acts freely and which forms a fiber bundle $E \twoheadrightarrow E/G$.
For this purpose we assume that $G$ is {\em linear},
i.e. that $G$ is a Lie subgroup of $GL(n)$, for some $n \in \BB N$.
Then we can make a trivial modification of the construction of the space $E$
given in Section 1.2 in \cite{GS}.

Let $L^2[0,\infty)$ denote the space of square integrable functions
on the positive real numbers relatively to the standard Lebesgue measure.
This space comes equipped with an inner product which gives it
the norm topology and makes it a Hilbert space.
Consider the space of $n$-frames
$$
E = \{ {\bf f} = (f_1, \dots,f_n) \in
L^2[0,\infty) \times \dots \times L^2[0,\infty)
;\: f_1, \dots, f_n \text{ are linearly independent}\}
$$
which is an open subset of $L^2[0,\infty) \times \dots \times L^2[0,\infty)$
equipped with the product topology.
We define the action of $G$ on $E$ by
$$
g \cdot (f_1, \dots,f_n) = (\tilde f_1, \dots, \tilde f_n),
\qquad
\tilde f_i = \sum_{j=1}^n a_{ij} f_j,
$$
where $g \in G \subset GL(n)$ is represented by an invertible
$n \times n$ matrix $(a_{ij})$.
Clearly, $G$ acts on $E$ continuously and freely.

Let $E' \subset E$ denote the subset of $n$-tuples of functions which all
vanish on the interval $[0,1]$.

\begin{lem}
The subset $E'$ is a deformation retract of $E$.
\end{lem}

\pf
For any $f \in L^2[0,\infty)$ define $T_tf$ by
$$
T_tf(x) =
\begin{cases}
0 & \text{for } 0 \le x < t;  \\
f(x-t) & \text{for } t \le x < \infty.
\end{cases}
$$
Define
$$
{\bf T}_t {\bf f} = (T_tf_1, \dots,T_tf_n) \qquad \text{for }
{\bf f} = (f_1, \dots,f_n) \in L^2[0,\infty) \times \dots \times L^2[0,\infty).
$$
Since ${\bf T}_t$ preserves linear independence we see that ${\bf T}_t$
is a deformation retract of $E$ into $E'$.
\qed

\begin{prop}
The space $E$ is contractible and the projection
$E \twoheadrightarrow E/G$ is a (locally trivial) fiber bundle.
\end{prop}

\pf
To prove that $E$ is contractible it is sufficient to show that
$E'$ is contractible to a point within $E$.
Pick an $n$-frame ${\bf g} = (g_1, \dots, g_n) \in E$ such that all its
components are supported in $[0,1]$.
For each ${\bf f} = (f_1, \dots , f_n) \in E'$ we define
$$
{\bf R}_t {\bf f} = \bigl( tg_1 + (1-t)f_1, \dots tg_n + (1-t)f_n \bigr),
\qquad t \in [0,1].
$$
Clearly, ${\bf R}_t (E') \subset E$ for all $t$
and ${\bf R}_t$ is a continuous deformation of $E'$ to ${\bf g}$ within $E$.

It remains to show that $p: E \twoheadrightarrow E/G$ is a fiber bundle.
Suppose first that $G=GL(n)$.
Pick a Hilbert space orthonormal basis $\{v_1,v_2,\dots\}$ of $L^2[0,\infty)$.
For each sequence of $n$ different integers
$1 \le i_1 < i_2 < \dots < i_n < \infty$ we set
$\BB R^n_{\{i_1,\dots,i_n\}} = \BB R$-span of $v_{i_1}, \dots, v_{i_n}$
in $L^2[0,\infty)$, and define
$$
U_{\{i_1,\dots,i_n\}} = \biggl\{ {\bf f} = (f_1, \dots,f_n) \in E ;\:
\begin{matrix}
\text{the orthogonal projections of $f_1,\dots,f_n$}  \\
\text{onto $\BB R^n_{\{i_1,\dots,i_n\}}$ form a basis in
$\BB R^n_{\{i_1,\dots,i_n\}}$}
\end{matrix} \biggr\}.
$$
The sets $U_{\{i_1,\dots,i_n\}}$ are $G$-invariant and form an open
covering of $E$. Hence their projections
$\{p(U_{\{i_1,\dots,i_n\}}) \}_{1 \le i_1 < i_2 < \dots < i_n < \infty}$
form an open covering of $E/G$.
Note that every element ${\bf f} \in U_{\{i_1,\dots,i_n\}}$ can be uniquely
written as ${\bf f} = g_{\bf f} \cdot {\bf f'}$, where $g_{\bf f} \in GL(n)$
and ${\bf f'} = (f_1',\dots,f_n') \in E$ is such that the orthogonal
projections of $f_1',\dots,f_n'$ onto $\BB R^n_{\{i_1,\dots,i_n\}}$ are 
$v_{i_1}, \dots, v_{i_n}$ respectively.
The maps
$\varphi_{\{i_1,\dots,i_n\}}: U_{\{i_1,\dots,i_n\}} \to
G \times p(U_{\{i_1,\dots,i_n\}})$ defined by
$$
\varphi_{\{i_1,\dots,i_n\}} : {\bf f} \mapsto (g_{\bf f}, p({\bf f})),
$$
provide a trivialization of $E \twoheadrightarrow E/G$ over each
$p(U_{\{i_1,\dots,i_n\}})$.

Finally, if $G$ is a proper subgroup of $GL(n)$, then 
$p: E \twoheadrightarrow E/G$ is a pullback of $E \twoheadrightarrow E/GL(n)$:
$$
\begin{matrix}
E & \longrightarrow & E  \\
\downarrow & \quad & \downarrow  \\
E/G & \longrightarrow & E/GL(n)
\end{matrix}
$$
hence a locally trivial fiber bundle too.
\qed

\separate

\begin{section}
{Restricting to subgroups}
\end{section}

Let $K \subset G$ be a Lie subgroup.
If $E$ is a contractible space on which $G$ acts freely and which forms a
fiber bundle $E \twoheadrightarrow E/G$,
then the pullback diagram of fiber bundles
$$
\begin{matrix}
E & \longrightarrow & E  \\
\downarrow & \quad & \downarrow  \\
E/K & \longrightarrow & E/G
\end{matrix}
$$
shows that $E$ ``works'' for $K$ as well.
Thus we have a map on equivariant cohomology
$$
res_K^G: \: H^*_G(X) \to  H^*_K(X)
$$
induced by the projection of spaces and the map on ordinary cohomology
$$
(X \times E)/K \twoheadrightarrow (X \times E)/G
\qquad \text{and} \qquad
H^*((X \times E)/G) \to H^*((X \times E)/K).
$$
This map on equivariant cohomology $res_K^G$ is canonical in the sense that
it does not depend on the choice of space $E$ (which can be proved by the
same argument that shows that Definition \ref{def} does not depend on the
choice of $E$).

Note that $(X \times E)/K \twoheadrightarrow (X \times E)/G$
is a fiber bundle with fiber over each point homeomorphic to
the homogeneous space $G/K$.
Hence we get a spectral sequence relating the equivariant cohomologies
$H^*_G(X)$, $H^*_K(X)$ and the ordinary cohomology $H^*(G/K)$.
In particular,

\begin{prop}
Suppose that the homogeneous space $G/K$ is contractible.
Then the restriction map
$$
res_K^G: \: H^*_G(X) \to  H^*_K(X)
$$
is an isomorphism.
\end{prop}

Two special cases are worth mentioning:

\begin{cor}
Suppose that the group $G$ is simply connected nilpotent, then
$$
H^*_G(X) \simeq  H^*(X).
$$
\end{cor}

\pf
Let $K$ be the trivial subgroup $\{e\}$. Since $G$ is simply connected
nilpotent, it is diffeomorphic to its Lie algebra via the exponential map
(see \cite{Kn}, for instance), hence contractible and
$$
res_{\{e\}}^G: \: H^*_G(X) \tilde \longrightarrow  H^*_{\{e\}}(X)
\simeq H^*(X).  \qquad \qquad \square
$$

\begin{cor}
Suppose that the group $G$ is linear reductive and $K \subset G$
is a maximal compact subgroup, then
$$
res_K^G: \: H^*_G(X) \to  H^*_K(X)
$$
is an isomorphism.
\end{cor}

\pf
Every reductive group has a Cartan decomposition (see \cite{Kn}, for instance)
which implies $G/K$ is contractible.
\qed

The last corollary implies that the twisted deRham complex
(see \cite{GS}, for instance) used to compute the $K$-equivariant cohomology 
can also be used to compute the $G$-equivariant cohomology.
It may appear at first that there is no interesting equivariant cohomology
theory for non-compact reductive groups since everything just gets reduced
to the action of their maximal compact subgroups.
We will see in Section \ref{new} that it is not so.

\separate

\begin{section}
{Equivariant forms and localization}
\end{section}

Equivariant forms were introduced in 1950 by Henri Cartan \cite{Ca1},
\cite{Ca2} (see also \cite{BGV} and \cite{GS}).
Let $G$ be a real Lie group acting on an oriented manifold $M$,
let $\g g$ be the Lie algebra of $K$, and let
$\Omega^*(M)$ denote the algebra of smooth differential forms on $M$.
Recall that a $G$-equivariant form is a map
$\alpha : \g g \to \Omega^*(M)$ such that
$$
\alpha(X) = g^{-1} \cdot \alpha (Ad(g) X)
\qquad \forall X \in \g g,\: g \in G.
$$
If the group $G$ is commutative (the circle group $S^1$ is a very
interesting example), then a $G$-equivariant form is really a $G$-invariant
map $\g g \to \Omega^*(M)$.

We define a {\em twisted deRham differential} by
$$
(d_{equiv} \alpha) (X) = d(\alpha(X)) + \iota(\VF_X) (\alpha(X)),
$$
where $d$ denotes the ordinary deRham differential and
$\iota(\VF_X)$ denotes contraction by the vector field $\VF_X$ induced by
the infinitesimal action of $X$ on $M$.
The map $d_{equiv}$ preserves $G$-equivariant forms and $(d_{equiv})^2=0$.

Equivariantly closed forms occur naturally in symplectic geometry.
If $M$ has a $G$-invariant symplectic form $\omega$ which admits a moment map
$\mu: M \to \g g^*$ for the action of $G$,
then $\omega + \mu$ and $e^{\omega+\mu}$ are equivariantly
closed forms. (See Chapter 9 of \cite{GS} for details.)

Let $\Omega_G^*(M)$ denote the complex of $G$-equivariant forms which
depend on $X \in \g g$ polynomially and the degree of an element 
$\alpha : \g g \to \Omega^*(M)$ is defined as twice its degree as a polynomial
on $\g g$ plus the degree of the differential form in $\Omega^*(M)$.
H.~Cartan proved in \cite{Ca1} and \cite{Ca2} (see also \cite{GS}) that,
when the group $G$ is compact, this complex $\Omega_G^*(M)$
computes the $G$-equivariant cohomology of $M$.
For this reason equivariant forms are very important and $\Omega_G^*(M)$
is called the twisted deRham complex.

\separate

If $\alpha$ is a possibly non-homogeneous equivariant differential form,
$\alpha_{[k]}$ denotes its homogeneous component of degree $k$.
For $X \in \g g$, we denote by $M^X$ the set of zeroes of the vector
field $\VF_X$ on $M$.
Let us assume for simplicity that the manifold $M$ is compact and the set
of zeroes $M^X$ is discrete (hence finite).
Then the integral localization formula can be stated as follows:

\begin{thm} {\rm (Theorem 7.11 in \cite{BGV})}  \label{BGV}
Suppose that the group $G$ is compact and
$\alpha : \g g \to \Omega^*(M)$ is an equivariant form
such that $d_{equiv} \alpha =0$ (i.e. $\alpha$ is equivariantly closed). Then
\begin{equation}  \label{BGVformula}
\int_M \alpha(X)_{[\dim M]} = (-2\pi)^{\dim M/2}
\sum_{p \in M^X} \frac {\alpha(X)_{[0]}(p)}{\det^{1/2} (L_p)},
\end{equation}
where $L_p: T_pM \to T_pM$ is a linear automorphism of the tangent space at $p$
induced by the Lie action of $-\VF_X$ on $M$, and $\det^{1/2} (L_p)$ is
a canonically defined (\cite{BGV}, Section 7.2) square root of the
determinant of this transformation depending only on the orientation of $T_pM$.
\end{thm}

This result relates a global object (i.e. integral) with
locally defined objects such as the quotients
$\frac {\alpha(X)_{[0]}(p)}{\det^{1/2} (L_p)}$ at the zeroes $p \in M^X$.
It was originally proved by N.~Berline and M.~Vergne \cite{BV1}
and independently by M.~Atiyah and R.~Bott \cite{AB}.
In symplectic geometry this result is often called the Duistermaat-Heckman
theorem.

There is a similar localization formula when the set of zeroes $M^X$
is not discrete, then the summation over $M^X$ is replaced with integration
\cite{BGV}. Compactness of $M$ is not essential either.
For example, the Fourier transform of a coadjoint orbit was originally
computed by W.~Rossmann \cite{Ro1}.
Then N.~Berline and M.~Vergne \cite{BV2} found an easier computation of
the Fourier transform essentially by applying (\ref{BGVformula}) to an
integral over the coadjoint orbit (which is not compact) and making sure
that everything decays fast at infinity and the integral converges.
This localization formula has many other applications;
we will mention one more in the next section.
Unfortunately, this formula fails when the acting Lie group $G$
is not compact: there simply may not be enough fixed points present.

\begin{ex}  {\em
Let $G = SL(2, \BB R)$. Then there exists an
essentially unique $SL(2, \BB R)$-invariant symmetric bilinear form $B$ on
$\g g = \g {sl}(2,\BB R)$, say $B(X,Y)$ is the Killing form
$\tr(ad(X)ad(Y))$.
This bilinear form $B$ induces an $SL(2,\BB R)$-equivariant isomorphism
$I: \g {sl}(2,\BB R) \, \tilde \to \, \g {sl}(2,\BB R)^*$.
Let ${\cal O} \subset \g {sl}(2,\BB R)^*$ denote the coadjoint orbit of
$I \begin{pmatrix} 0 & 1 \\ -1 & 0 \end{pmatrix}$.
Like all coadjoint orbits, ${\cal O}$ possesses a canonical
symplectic form $\omega$ which is the top degree part of a certain
equivariantly closed 2-form.
Although ${\cal O}$ is not compact, the symplectic volume
$\int_{\cal O} \omega$ still exists as a distribution on $\g {sl}(2,\BB R)$.
Let $\g {sl} (2,\BB R)'_{split} \subset \g {sl}(2,\BB R)$ be the
open subset consisting of $X \in \g {sl}(2,\BB R)$
with distinct real eigenvalues.
Now, if we take any element $X \in \g {sl} (2,\BB R)'_{split}$,
then one can see that the vector field on ${\cal O}$ generated by $X$
has no zeroes, i.e. ${\cal O}^X = \varnothing$.
Thus if there were a fixed point integral localization
formula like in the case of compact groups,
this formula would suggest that the distribution
determined by $\int_{\cal O} \omega$ vanishes on the open set
$\g {sl} (2,\BB R)'_{split}$.
But it is known that the restriction $\int_{\cal O} \omega$ to
$\g {sl} (2,\BB R)'_{split}$ is {\em not} zero.
\qed
}\end{ex}

On the other hand, recent results from representation theory, namely
the two character formulas for representations of reductive Lie groups
due to M.~Kashiwara, W.~Rossmann, W.~Schmid and K.~Vilonen described in
\cite{Sch}, \cite{SchV2} strongly suggest that integral localization
formula should extend to actions of non-compact groups.
The above example illustrates some of the obvious obstacles to having
a localization formula when the acting group $G$ is not compact:
\begin{itemize}
\item
In order to have a truly new result where the action of $G$ does not
factor through action of some compact group we must allow non-compact
manifolds or Borel-Moore homology cycles with infinite support.
But then we need to worry about convergence of the integral.
We will resolve this problem by restricting the class of forms that we will
integrate and by introducing a new (weaker) notion of convergence of
integrals in the sense of distributions on $\g g$.
\item
For an arbitrary cycle with infinite support, the fixed points tend to
``run away to infinity.'' This happens in total analogy with the failure
of the Lefschetz fixed point formula for non-compact spaces.
We will describe a class of cycles for which all fixed points are accounted
for. Even if a cycle does not contain enough fixed points, it may be possible
to deform it into a new cycle for which the localization formula is known
to be true and the integral stays unchanged. I make an attempt to study such
deformations in \cite{L3}.
\end{itemize}

\separate

\begin{section}
{Some representation theory}
\end{section}

Recall that $\g g$ is the Lie algebra of a Lie group $G$.
Let $\g g_{\BB C} = \g g \otimes_{\BB R} \BB C$
be the complexified Lie algebra. Then the flag variety ${\cal B}$
of $\g g_{\BB C}$ is the variety of all Borel (maximal solvable)
subalgebras of $\g g_{\BB C}$. 
The space ${\cal B}$ is a smooth complex projective variety.
If the group $G$ is compact, then all irreducible representations
of $G$ can be enumerated by their {\em highest weights} $\lambda$
lying in the {\em weight lattice} $\Lambda \subset i \g g^*$.

The Borel-Weil-Bott theorem \cite{Bott} can be regarded as an explicit
construction of a holomorphic $G$-equivariant line bundle
${\cal L}_{\lambda} \to {\cal B}$ such that the resulting
representation of $G$ in the cohomology groups is:
\begin{eqnarray*}
&& H^p({\cal B}, {\cal O}({\cal L}_{\lambda}))=0 \text{\quad if $p \ne 0$}, \\
&& H^0({\cal B}, {\cal O}({\cal L}_{\lambda})) \simeq \pi_{\lambda},
\end{eqnarray*}
where $\pi_{\lambda}$ denotes the irreducible representation of
$G$ of highest weight $\lambda$.
Then N.~Berline and M.~Vergne \cite{BV1}, \cite{BGV} observed that the
character of $\pi_{\lambda}$, as a function on $\g g$, can be expressed as
an integral over ${\cal B}$ of a certain naturally defined equivariantly
closed form. They proved it by applying their localization formula
(\ref{BGVformula}) and matching contributions from fixed points with terms
of the Weyl character formula.
(This is a restatement of Kirillov's character formula.)

\separate

Next we describe some recent results on representations of non-compact groups.
We fix a connected complex algebraic linear reductive Lie group $G_{\BB C}$
which is defined over $\BB R$.
We will be primarily interested in a real Lie subgroup
$G \subset G_{\BB C}$ lying between the group of real points $G_{\BB C}(\BB R)$
and the identity component $G_{\BB C}(\BB R)^0$.
We regard $G$ as a real reductive Lie group
(e.g. $SL(n,\BB R)$, $GL(n, \BB R)$, $U(n)$, $Sp(n)$,\dots).
Because there may not be enough finite-dimensional representations,
we consider topological vector spaces $V$ of possibly infinite dimension
with continuous $G$-action.
A reasonable category of representations consists of
{\em admissible} representations of {\em finite length}.
(A representation $\pi$ has finite length if every increasing chain
of closed, invariant subspaces breaks off after finitely many steps;
$\pi$ is admissible if its restriction to a maximal compact subgroup $K$
contains any irreducible representation of $K$ at most finitely often.)
An irreducible unitary representation is always of this kind.
Although trace of a linear operator in an infinite-dimensional
space cannot be defined in general, it is still possible to
define a character $\theta_{\pi}$ as an
$Ad(G)$-invariant distribution on $\g g$.

M.~Kashiwara and W.~Schmid \cite{KaSchm}
generalize the Borel-Weil-Bott construction.
Instead of line bundles over the flag variety ${\cal B}$ they consider
$G$-equivariant sheaves ${\cal F}$ and, for each integer
$p \in \BB Z$, they define representations of $G$ in
$\operatorname{Ext}^p({\cal F},{\cal O}_{\cal B})$.
Such representations turn out to be admissible of finite length.
Let $\theta$ be the character of the virtual representation of $G$
$$
\sum_p (-1)^p \operatorname{Ext}^p
(\BB D {\cal F}, {\cal O}_{\cal B}(\lambda)),
$$
where $\BB D {\cal F}$ denotes the Verdier dual of ${\cal F}$
and $\lambda$ is some twisting parameter lying in $\g h_{\BB C}^*$ --
the dual space of the universal Cartan algebra of $\g g_{\BB C}$.

Then W.~Schmid and K.~Vilonen \cite{SchV2} prove
two character formulas for $\theta$.
Since $\theta$ is a distribution on $\g g$,
let $\varphi \in {\cal C}_c^{\infty}(\g g)$ be a test function on $\g g$
and let $dX$ be the Lebesgue measure on $\g g$.
We also let $T^*{\cal B}$ denote the cotangent space of ${\cal B}$.
Then the integral character formula says
$$
\theta(\varphi) = \int_{Ch({\cal F})} \alpha(\varphi),
$$
where $\alpha(\varphi)$ is an equivariantly closed form on $T^*{\cal B}$
(which does not depend on ${\cal F}$), and
$Ch({\cal F})$ is the characteristic cycle of ${\cal F}$ which lies in
$T^*{\cal B}$.
(Characteristic cycles of constructible sheaves were introduced by
M.~Kashiwara and their definition can be found in
\cite{KaScha}, \cite{SchV1}, \cite{Schu}.)
On the other hand, the fixed point character formula says
$$
\theta(\varphi) = \int_{\g g} F_{\theta} \varphi \, dX,
$$
where $F_{\theta}$ is an $Ad(G)$-invariant locally $L^1$-function
such that its restriction to the set of regular semisimple elements of
$\g g$ can be represented by an analytic function.
The value of this analytic function at a regular semisimple element
$X \in \g g$ is given by the formula
$$
F_{\theta}(X) = (-2\pi)^{\dim_{\BB R} {\cal B}/2}
\sum_{p \in {\cal B}^X} m_p(X) \cdot Q(X),
$$
where $Q(X)$ is the term which appeared both on the right hand side of
(\ref{BGVformula}) and in the Weyl character formula,
and $m_p(X)$ is a certain integer multiplicity which
is exactly the local contribution of $p$ to the Lefschetz
fixed point formula, as generalized to sheaf cohomology by
M.~Goresky and R.~MacPherson \cite{GM}.
These multiplicities are determined in terms of
local cohomology of ${\cal F}$.

In the special case when the group $G$ is compact, the former reduces to
Kirillov's character formula and the latter -- the Weyl character formula.
The fixed point formula was conjectured by M.~Kashiwara
\cite{Ka}, and its proof uses a generalization of the Lefschetz fixed
point formula to sheaf cohomology \cite{GM}.
On the other hand, W.~Rossmann \cite{Ro1} established existence
of an integral character formula over an unspecified
Borel-Moore cycle.
These character formulas are proved in \cite{SchV2} independently of each
other using representation theory methods.

\separate

\begin{section}
{New localization formula}  \label{new}
\end{section}

As in the previous section, we fix a connected complex algebraic linear
reductive Lie group $G_{\BB C}$ which is defined over $\BB R$, let
$G$ be a real Lie subgroup of $G_{\BB C}$
lying between the group of real points $G_{\BB C}(\BB R)$
and the identity component $G_{\BB C}(\BB R)^0$,
and regard $G$ as a real reductive Lie group.
Our ambient space will be the holomorphic cotangent space $T^*M$ of a smooth
complex projective variety $M$ on which $G_{\BB C}$ acts algebraically.
We will also assume that any maximal complex torus
$T_{\BB C} \subset G_{\BB C}$ acts on $M$ with isolated fixed points.
Then there are only finitely many of them because $M$ is compact.
(This condition is satisfied in all applications we have in mind.)

Let $\sigma$ denote the canonical complex algebraic holomorphic symplectic
form on $T^*M$.
The Borel-Moore cycles $\Lambda \subset T^*M$ over which we will integrate
will be subject to the following three properties:
\begin{itemize}
\item
$\Lambda$ is $G$-invariant;
\item
$\Lambda$ is real Lagrangian, i.e. $\re \sigma |_{\Lambda} \equiv 0$
and $\dim_{\BB R} \Lambda = \dim_{\BB R} M$;
\item
$\Lambda$ is {\em conic}, i.e. invariant under the scaling action of
positive reals $\BB R^{>0}$ on $T^*M$ (but not necessarily under
the actions of $\BB C^{\times}$ or $\BB R^{\times}$).
\end{itemize}

\begin{ex}
Let $N \subset M$ be a closed $G$-invariant real
submanifold, and let $\Lambda$ be the conormal space $T^*_NM$ equipped
with some orientation.

An interesting example is
$G = GL(n,\BB R) \subset GL(n, \BB C) = G_{\BB C}$
acting naturally on a complex Grassmanian $Gr_{\BB C}(k,n)$.
Let $N$ be the real Grassmanian $Gr_{\BB R}(k,n)$ sitting inside
$Gr_{\BB C}(k,n)$ and $\Lambda = T^*_{Gr_{\BB R}(k,n)} Gr_{\BB C}(k,n)$.
\end{ex}

\begin{rem}
Any such cycle $\Lambda$ can be realized as a characteristic cycle
$Ch({\cal F})$ of some $G$-equivariant constructible sheaf ${\cal F}$
(\cite{KaScha}, \cite{SchV1}, \cite{Schu}).
\end{rem}

Let $U$ be another subgroup of $G_{\BB C}$ such that,
letting $\g u$ denote the Lie algebra of $U$, we have an isomorphism
$\g u \otimes_{\BB R} \BB C \simeq \g g_{\BB C}$.
For instance, $U$ may equal $G$, but in all applications we have in mind
$U$ is a compact real form of $G_{\BB C}$, i.e. a maximal compact subgroup.
We denote by $\Omega^{(p,q)}(M)$ the space of complex-valued
differential forms of type $(p,q)$ on $M$.
We will consider forms $\alpha : \g g_{\BB C} \to \Omega^*(M)$
which satisfy the following three conditions:
\begin{itemize}
\item
The assignment $X \mapsto \alpha(X) \in \Omega^*(M)$
depends holomorphically on $X \in \g g_{\BB C}$;
\item
For each $k \in \BB N$ and each $X \in \g g_{\BB C}$,
$$
\alpha(X)_{[2k]} \in
\bigoplus_{\begin{matrix} p+q=2k \\ p \ge q \end{matrix}} \Omega^{(p,q)}(M);
$$
\item
The restriction of $\alpha$ to $\g u \subset \g g_{\BB C}$
is an equivariantly closed form with respect to $U$.
\end{itemize}

\begin{ex}  {\em
A $U$-equivariant characteristic form $\alpha: \g u \to \Omega^*(M)$
associated to a $U$-equivariant vector bundle over $M$
(see Section 7.1 of \cite{BGV}) satisfies the third condition.
Since it depends on $X \in \g u$ polynomially,
$\alpha$ extends uniquely to a map $\alpha: \g g_{\BB C} \to \Omega^*(M)$
so that the first condition is satisfied.
Finally, for each $X \in \g g_{\BB C}$,
$$
\alpha(X) \in \bigoplus_k \Omega^{(k,k)}(M),
$$
so that the second condition is satisfied too.
This is the most important class of forms satisfying these conditions.
\qed
}\end{ex}

Let $\mu: T^*M \to \g g_{\BB C}^*$ be the ordinary moment map:
$$
\mu(\xi): X \mapsto \langle \xi, \VF_X \rangle,
\qquad \xi \in T^*M, \: X \in \g g_{\BB C}.
$$
The integrals will be defined as distributions on $\g g$,
so let $\varphi \in {\cal C}_c^{\infty}(\g g)$ be a test function,
and let $dX$ denote the Lebesgue measure on $\g g$.
The new localization formula will apply to integrals
of the following kind:
\begin{equation}  \label{int}
\int_{\Lambda} \Bigl( \int_{\g g} 
e^{\langle X, \mu(\xi) \rangle + \sigma}
\wedge \varphi(X) \alpha(X) \,dX \Bigr)_{[\dim_{\BB R}M]},
\qquad X \in \g g ,\: \xi \in |\Lambda| \subset T^*M.
\end{equation}
The idea to consider integrals of this kind was inspired by the shape
of the integral character formula described in the previous section.
The inside integral
$$
\int_{\g g} e^{\langle X, \mu(\xi) \rangle + \sigma}
\wedge \varphi(X) \alpha(X) \,dX
$$
is essentially the Fourier transform of $\varphi(X) \alpha(X)$ which decays
rapidly in the imaginary directions of
$\g g_{\BB C}^* \simeq \g g^* \oplus i \g g^*$.
We denote by
$$
\supp(\sigma|_{\Lambda})
$$
the closure in $T^*M$ of the set of smooth points of the support
$|\Lambda|$ where $\sigma|_{|\Lambda|} \ne 0$.
Then integral (\ref{int}) converges provided that the moment map $\mu$
is proper on $\supp(\sigma|_{\Lambda})$.
In particular, (\ref{int}) is well-defined when $\mu$ is proper on $|\Lambda|$.

Now the main result of \cite{L4} says that if the support of $\varphi$
lies in $\g g'$ ($\g g$ without a finite number of certain hypersurfaces)
then integral (\ref{int}) can be rewritten as
$$
\int_{\Lambda} \Bigl( \int_{\g g} 
e^{\langle X, \mu(\xi) \rangle + \sigma}
\wedge \varphi(X) \alpha(X) \,dX \Bigr)_{[\dim_{\BB R}M]}
= \int_{\g g} F_{\alpha}(X) \varphi(X) \,dX,
$$
where $F_{\alpha}$ is an $Ad(G \cap U)$-invariant function on $\g g'$
given by the formula
\begin{equation}  \label{mainequation}
F_{\alpha}(X) = (-2\pi)^{\dim_{\BB R} M/2}
\sum_{p \in M^X} m_p(X) \frac {\alpha(X)_{[0]}(p)}{\det^{1/2} (L_p)}.
\end{equation}
As in the fixed point character formula,
$m_p(X)$ is an integer multiplicity equal the local contribution of $p$
to the Lefschetz fixed point formula, as generalized to sheaf cohomology by
M.~Goresky and R.~MacPherson \cite{GM}.
These multiplicities are determined in \cite{L4} in terms of
local cohomology of ${\cal F}$, where ${\cal F}$ is any sheaf
with characteristic cycle $Ch({\cal F}) = \Lambda$.
The set $\g g'$ is essentially the set of regular semisimple elements
of $\g g$ on which the denominators $\det^{1/2} (L_p)$ do not vanish.

\begin{rem}
In the special case when $\Lambda = M$ as oriented cycles,
$\Lambda$ is $U$-invariant, each multiplicity
$m_p(X)$ equals 1 and this theorem can be easily deduced from
the classical integral localization formula (Theorem \ref{BGV}).

Notice that the cycle $\Lambda$ is invariant with respect to the
action of the group $G$ which need not be compact, while the form 
$\alpha: \g g_{\BB C} \to \Omega^*(M)$ is required to be equivariant
with respect to a different group $U$ only, and $U$
may not preserve the cycle $\Lambda$.
\end{rem}

There are several important applications of
the new localization formula (\ref{mainequation}).
In \cite{L1} I give a geometric proof of the integral character formula
by matching the contributions of the fixed points with the terms of
the fixed point character formula.
Article \cite{L2} gives a very accessible introduction to \cite{L1}
and explains the key ideas used there by way of examples and illustrations.
Another application of (\ref{mainequation}) is a generalization of
the Gauss-Bonnet theorem to sheaf cohomology \cite{L4}:
the Euler characteristic $\chi(M,{\cal F})$ can be expressed as an
integral over the characteristic cycle $Ch({\cal F})$ of an extension of
the equivariant Euler form.
This formula is proved by comparing the contributions of the fixed points
with M.~Kashiwara's generalization of the Hopf index theorem stated as
Corollary 9.5.2 in \cite{KaScha},
$$
\chi(M, {\cal F}) = \# \bigl( [M] \cap Ch({\cal F}) \bigr),
$$
where $[M]$ is the fundamental cycle of $M$.
In \cite{L5} I use (\ref{mainequation}) to prove a Riemann-Roch-Hirzebruch
type integral formula for characters of representations of reductive groups.
These results strongly suggest that many statements which previously
were known in the compact group setting only can be generalized to
non-compact groups.
\separate

The proof of (\ref{mainequation}) utilizes a combination of two deformations
and its idea can be outlined as follows:
\begin{itemize}
\item
The integrand of (\ref{int}) is a closed differential form (easy).

\item
We introduce the first deformation
$$
\Theta_t(X): T^*M \to T^*M,
\qquad \qquad X \in \g g, \: t \in [0,1],
$$
for the purpose of making the integral
\begin{equation}  \label{deformedint}
\int_{\g g \times \Theta_t(X)_* (\Lambda)}
\Bigl( e^{\langle X, \mu(\xi) \rangle + \sigma}
\wedge \varphi(X) \alpha(X) \Bigr)_{[\dim_{\BB R}M]} \,dX
\end{equation}
absolutely convergent for $t \in (0,1]$, and $\Theta_0(X)$ is the identity map.

\item
The crux of the matter is that the integrals (\ref{int}) and
(\ref{deformedint}) are equal, i.e. the integral (\ref{int}) stays unchanged
after this deformation.
This statement is not at all obvious since the integral (\ref{int})
is not absolutely convergent.

\item
For each sufficiently regular $X \in \g g$,
the cycle $\Lambda$ is homologous to
$$
\sum_{p \in M^X} m_p(X) \cdot T^*_pM
$$
inside the set
\begin{equation}  \label{set}
\{ \xi \in T^*M ;\: \re \langle X, \mu(\xi) \rangle \le 0 \}
\end{equation}
(linearization theorem).
Here each cotangent space $T^*_pM$ is given some orientation.
Notice that the expression
$\langle X, \mu(\xi) \rangle$ appears in
the exponent of the integrand, thus one should expect good behavior of the
integral as the cycle is deformed inside the set (\ref{set}).
Combining this with the first deformation, we obtain a deformation of
$\Theta_t(X)_* (\Lambda)$ into
$$
\sum_{p \in M^X} m_p(X) \cdot \Theta_t(X)_* (T^*_pM).
$$

\item
The integral (\ref{deformedint}) stays unchanged during the above deformation
of $\Theta_t(X)_* (\Lambda)$. This is essentially because the integrand is
a closed form and the convergence of (\ref{deformedint}) is absolute.

\item
The contribution of each cycle $\Theta_t(X)_* (T^*_pM)$ to the integral
is exactly 
$$
(-2\pi)^{\dim_{\BB R} M/2} \int_{\g g}
\frac {\alpha(X)_{[0]}(p)}{\det^{1/2} (L_p)} \varphi(X) \,dX.
$$
\end{itemize}

\begin{section}
{Duistermaat-Heckman measures}
\end{section}

As before, $G$ is a linear real reductive Lie group with complexification
$G_{\BB C}$, we denote by $\g g$ and $\g g_{\BB C}$ their respective
Lie algebras.
We pick another subgroup $U$ of $G_{\BB C}$ such that, letting
$\g u$ be the Lie algebra of $U$, we have an isomorphism
$\g u \otimes_{\BB R} \BB C \simeq \g g_{\BB C}$.
For instance, $U$ may equal $G$, but in most interesting
situations $U$ is a compact real form of $G_{\BB C}$
(i.e. a maximal compact subgroup).

Let $M$ be a smooth complex projective variety equipped with an algebraic
action of $G_{\BB C}$ preserving a complex-valued 2-form $\omega$,
and suppose that the restriction of the
$G_{\BB C}$-action to $U$ is Hamiltonian with respect to $\omega$.
In other words, there exists a moment map
$\phi: M \to \g u^* \otimes_{\BB R} \BB C \simeq \g g_{\BB C}^*$
such that
$$
\iota(\operatorname{VF}_X) \omega = -d\phi(X),
\qquad \forall X \in \g u,
$$
where $\iota(\VF_X)$ denotes contraction by the vector field $\VF_X$
induced by the infinitesimal action of $X$ on $M$.
Note that we do not require the 2-form $\omega$ to be symplectic,
i.e. $\omega^{\dim_{\BB R} M/2} \ne 0$.
Even the case $\omega=0$, $\phi=0$ is interesting enough,
but, of course, symplectic forms are the most interesting ones.
We can regard $\phi: M \to \g g_{\BB C}^*$ as a map
$\phi: \g g_{\BB C} \to {\cal C}^{\infty}(M)$.
Then $\omega +\phi$ is an equivariantly closed form on $M$
for the action of $U$.

Recall that $\sigma$ denotes the canonical complex algebraic holomorphic
symplectic form on the holomorphic cotangent bundle $T^*M$, and
$\mu: T^*M \to \g g_{\BB C}^*$ is the ordinary moment map.
As in the previous section, we fix a $G$-invariant real-Lagrangian cycle
$\Lambda \subset T^*M$ which is conic with respect to the scaling action
of $\BB R^{>0}$.

Set $n = \dim_{\BB C} M$. The Liouville form
$$
\frac {(\omega+\sigma)^n}{n!} =
\bigl( \exp (\omega+\sigma) \bigr)_{[2n]}
$$
determines a measure $\beta_{\Lambda}$ on $\Lambda$.
We will call the pushforward of this measure
$(\phi+\mu)_*(\beta_{\Lambda})$ on $\g g_{\BB C}^*$
the Duistermaat-Heckman measure.
That is, for a compactly supported smooth function
$f \in {\cal C}^{\infty}_c (\g g_{\BB C}^*)$,
\begin{equation}  \label{measure}
\int_{\g g_{\BB C}^*} f \, d(\phi+\mu)_*(\beta_{\Lambda}) \quad =_{def}
\quad \int_{\Lambda} \frac {(\omega+\sigma)^n}{n!} f \circ (\phi+\mu).
\end{equation}
The right hand side of (\ref{measure}) converges whenever the map
$\phi+\mu$ is proper on the set $\supp(\sigma|_{\Lambda})$.
This happens whenever $\mu$ is proper on $\supp(\sigma|_{\Lambda})$.
In particular, the pushforward $(\phi+\mu)_*(\beta_{\Lambda})$
is well-defined when $\mu$ is proper on $|\Lambda|$.

Duistermaat-Heckman measures are important invariants of symplectic manifolds
and there are so many papers on this subject that it is impossible to list
them all.
At first an explicit formula was given by J.~J.~Duistermaat and G.~J.~Heckman
\cite{DH} using the method of exact stationary phase in the special case when
$G$ is a compact torus acting with isolated fixed points.
It was extended to compact non-abelian groups by V.~Guillemin and E.~Prato
\cite{GP}. Then it was extended to compact non-abelian groups acting with
possibly non-isolated fixed points by L.~Jeffrey and F.~Kirwan \cite{JK}.
Many recent results on Duistermaat-Heckman measures are obtained by
computing their Fourier transforms using the integral localization formula
and then inverting these Fourier transforms.

Since the cycle $\Lambda$ is $G$-invariant, the moment map
$\mu$ takes purely imaginary values on its support $|\Lambda|$:
$$
\mu(|\Lambda|) \quad \subset \quad i\g g^* \quad \subset
\quad \g g^* \oplus i \g g^* \quad \simeq \quad
\g g_{\BB C}^*.
$$
Since $M$ is compact, the support of $(\phi+\mu)_*(\beta_{\Lambda})$,
which must be a subset of $(\phi+\mu)_*(|\Lambda|)$, is a subset of 
$\g g_{\BB C}^* \simeq \g g^* \oplus i \g g^*$
with bounded real part.

The Fourier transform of the Duistermaat-Heckman measure is a distribution
on $\g g$, i.e. a continuous linear functional on the space of
test functions ${\cal C}_c^{\infty}(\g g)$.
Following the conventions of \cite{L1}, \cite{L4} and \cite{SchV2} we define
the Fourier transform of $\varphi \in {\cal C}_c^{\infty}(\g g)$ without
the customary factor of $i = \sqrt{-1}$ in the exponent:
$$
\hat \varphi(\zeta) =
\int_{\g g} e^{\langle X, \zeta \rangle} \varphi(X)\, dX,
\qquad X \in \g g, \: \zeta \in \g g_{\BB C}^*.
$$
Notice that $\hat \varphi(\zeta)$ decays rapidly as $\zeta \to \infty$ and
the real part of $\zeta$ stays uniformly bounded.
Hence the value of the Fourier transform of $(\phi+\mu)_*(\beta_{\Lambda})$
at $\varphi \in {\cal C}_c^{\infty}(\g g)$ is
\begin{multline}  \label{FT1}
\widehat{(\phi+\mu)_*(\beta_{\Lambda})}(\varphi) =
\int \Bigl( \int_{\g g} e^{\langle X, \xi \rangle} \varphi(X) \,dX \Bigr)
\, d(\phi+\mu)_*(\beta_{\Lambda})   \\
= \int_{\Lambda} \Bigl( \int_{\g g}
e^{\langle X, \phi + \mu(\xi) \rangle} \varphi(X) \,dX \Bigr)
\frac {(\omega+\sigma)^n}{n!},
\qquad X \in \g g, \: \xi \in |\Lambda| \subset T^*M.
\end{multline}
We introduce a $U$-equivariant form
$\alpha: \g g_{\BB C} \to \Omega^*(M)$:
$$
\alpha(X) = \exp(\phi(X)+ \omega),
$$
then formula (\ref{FT1}) can be rewritten as
\begin{equation}  \label{FT2}
\widehat{(\phi+\mu)_*(\beta_{\Lambda})}(\varphi) =
\int_{\Lambda} \Bigl( \int_{\g g} 
e^{\langle X, \mu(\xi) \rangle + \sigma}
\wedge \varphi(X) \alpha(X) \,dX \Bigr)_{[\dim_{\BB R}M]},
\qquad X \in \g g ,\: \xi \in |\Lambda| \subset T^*M.
\end{equation}
This integral is exactly of type (\ref{int}), hence convergent
(essentially because $\hat \varphi (\zeta)$ decays rapidly as
$\zeta \to \infty$ and $\zeta \in (\phi+\mu)_*(|\Lambda|)$).
Since we want to apply the generalized localization formula
(\ref{mainequation}), we assume that any maximal complex torus
$T_{\BB C} \subset G_{\BB C}$ acts on $M$ with finitely many isolated
fixed points and that
$$
\omega \in \Omega^{(2,0)}(M) \oplus \Omega^{(1,1)}(M).
$$
Then (\ref{mainequation}) says that the restriction of the Fourier
transform of the Duistermaat-Heckman measure (\ref{FT2}) to
$\g g'$ (a certain open dense subset of $\g g$) equals
$$
\widehat{(\phi+\mu)_*(\beta_{\Lambda})}(\varphi) =
\int_{\g g_{\BB R}} F_{\omega}(X) \varphi(X) \,dX,
$$
where $F_{\omega}$ is an $Ad(G \cap U)$-invariant
function on $\g g'$ given by the formula
$$
F_{\omega}(X) = (-2\pi)^{\dim_{\BB R} M/2}
\sum_{p \in M^X} m_p(X) \frac {e^{\langle X, \phi(p) \rangle}}
{\det^{1/2} (L_p)}.
$$
where $M^X$ is the set of zeroes of the vector field $\operatorname{VF}_X$
on $M$, and $m_p(X)$'s are certain integer multiplicities.

\separate

\separate

\noindent
{\em E-mail address:} {matvei@math.umass.edu}

\noindent
{\em Department of Mathematics and Statistics, University of Massachusetts,
Lederle Graduate Research Tower, 710 North Pleasant Street, Amherst,
MA 01003}


\begin{thebibliography}{[KaSchm]}
\bibitem[AB]{AB} M.~F.~Atiyah and R.~Bott,
{\em The moment map and equivariant cohomology},
Topology {\bf 23} (1984), no. 1, 1-28.
\bibitem[BGV]{BGV} N. Berline, E. Getzler, M. Vergne,
{\em Heat Kernels and Dirac Operators}, Springer-Verlag, 1992.
\bibitem[BV1]{BV1} N.~Berline and M.~Vergne,
{\em Classes caract\'eristiques \'equivariantes.
Formules de localisation en cohomologie \'equivariante},
C. R. Acad. Sci. Paris {\bf 295} (1982), 539-541.
\bibitem[BV2]{BV2} N. Berline, M. Vergne,
{\em Fourier Transform of Orbits of the Coadjoint Representation},
Representation Theory of Reductive Lie Groups,
Progress in Mathematics, vol 40, Birkh\"auser, Basel, 1983, 53-57.
\bibitem[Bott]{Bott} R.~Bott, {\em Homogeneous vector bundles},
Annals of Math. {\bf 66} (1957), 203-248.
\bibitem[Ca1]{Ca1} H.~Cartan, {\em Notions d'alg\`ebre diff\'erentielle;
application aux groupes de Lie et aux vari\'et\'es o\`u op\`ere un groupe
de Lie}, Colloque de Topologie, C.B.R.M., Bruxelles 15-27 (1950).
\bibitem[Ca2]{Ca2} H.~Cartan, {\em La transgression dans un groupe de Lie
et dans un espace fibr\'e principal}.
Colloque de Topologie, C.B.R.M., Bruxelles 57-71 (1950).
\bibitem[DH]{DH} J.~J.~Duistermaat and G.~J.~Heckman,
{\em On the variation in the cohomology of the symplectic form
of the reduced phase space}, Invent. Math. {\bf 69} (1982), 259-268;
{\em Addendum}, ibid. {\bf 72} (1983), 153-158.
\bibitem[GM]{GM} M.~Goresky and R.~MacPherson,
{\em Local contribution to the Lefschetz fixed point formula},
Inventiones Math. {\bf 111} (1993), 1-33.
\bibitem[GP]{GP} V.~Guillemin and E.~Prato,
{\em Heckman, Kostant, and Steinberg Formulas for Symplectic Manifolds},
Advances in Math, {\bf 82} (1990), 160-179.
\bibitem[GS]{GS} V. Guillemin and S. Sternberg,
{\em Supersymmetry and Equivariant de Rham Theory},
Springer-Verlag, 1999.
\bibitem[JK]{JK} L.~Jeffrey and F.~Kirwan,
{\em Localization for non-abelian group actions},
Topology {\bf 34} (1995), 291-327.
\bibitem[Ka]{Ka} M.~Kashiwara, {\em Character, character cycle,
fixed point theorem, and group representations},
Advanced Studies in Pure Mathematics, vol. 14, Kinokuniya, Tokyo,
1988, 369-378.
\bibitem[KaScha]{KaScha} M.~Kashiwara and P.~Schapira,
{\em Sheaves on Manifolds}, Springer, 1990.
\bibitem[KaSchm]{KaSchm} M.~Kashiwara and W.~Schmid,
{\em Quasi-equivariant ${\cal D}$-modules, equivariant derived category, and
representations of reductive Lie groups},
Lie Theory and Geometry, in Honor of Bertram Kostant,
Progress in Mathematics, vol. 123, Birkh\"auser, Boston, 1994, pp. 457-488.
\bibitem[Kn]{Kn} A.~Knapp, {\em Lie Groups Beyond an Introduction},
Progress in Mathematics, vol. 140, Birkh\"auser, 2002.
\bibitem[L1]{L1} M.~Libine,
{\em A localization argument for characters of reductive Lie groups},
Jour. Func. Anal. {\bf 203} (2003), 197-236; also math.RT/0206019.
\bibitem[L2]{L2} M.~Libine,
{\em A Localization argument for characters of reductive Lie groups:
an introduction and examples},
Non-Commutative Harmonic Analysis: In Honor of Jacques Carmona,
Progress in Mathematics, vol. 220, Birkh\"auser, 2004, pp. 375-394;
also math.RT/0208024.
\bibitem[L3]{L3} M.~Libine,
{\em A fixed point localization formula for the
Fourier transform of regular semisimple coadjoint orbits},
math.DG/0302352, 2003, to appear in the Jour. Func. Anal.
\bibitem[L4]{L4} M.~Libine,
{\em Integrals of equivariant forms and
a Gauss-Bonnet theorem for constructible sheaves},
math.DG/0306152, 2003.
\bibitem[L5]{L5} M.~Libine, {\em Riemann-Roch-Hirzebruch integral formula
for characters of reductive Lie groups}, math.RT/0312454, 2003.
\bibitem[LSch]{LSch} M.~Libine and W.~Schmid,
{\em Geometric methods in representation theory}, math.RT/0402306, 2003,
to appear in the Proceedings of the International Euroschool and
Euroconference PQR2003.
\bibitem[Ro1]{Ro1} W.~Rossmann, {\em Kirillov's Character Formula for Reductive
Lie Groups}, Invent. Math. {\bf 48} (1978), 207-220.
\bibitem[Ro2]{Ro2} W. Rossmann, {\em Invariant Eigendistributions on a
Semisimple Lie Algebra and Homology Classes on the Conormal Variety I, II},
Jour. Func. Anal. {\bf 96} (1991), 130-193.
\bibitem[Sch]{Sch} W.~Schmid,
{\em Character formulas and localization of integrals},
Deformation Theory and Symplectic Geometry, Mathematical Physics Studies,
{\bf 20} (1997), Kluwer Academic Publishers, 259-270.
\bibitem[SchV1]{SchV1} W. Schmid and K. Vilonen,
{\em Characteristic cycles of constructible sheaves},
Inventiones Math. {\bf 124} (1996), 451-502.
\bibitem[SchV2]{SchV2} W. Schmid and K. Vilonen,
{\em Two geometric character formulas for reductive Lie groups},
Jour. AMS {\bf 11} (1998), 799-876.
\bibitem[{Sch\"u}]{Schu} J.~Sch\"urmann,
{\em Topology of Singular Spaces and Constructible Sheaves},
Monografie Matematyczne, vol. 63, Birkh\"auser, 2003.
\end{thebibliography}
\end{document}